\documentclass{elsart}

\usepackage{amsmath}
\usepackage{amsfonts}
\usepackage{amssymb}
\usepackage{latexsym,makeidx}
\usepackage{amsfonts}
\usepackage{amssymb}
\usepackage{epsfig}

\def\be{\begin{equation}}
\def\ee{\end{equation}}
\def\bc{\begin{center}}
\def\ec{\end{center}}

\newtheorem{lm}{Lemma}
\newtheorem{theorem}{Theorem}
\newtheorem{definition}{Definition}

\begin{document}

\begin{frontmatter}

\title{Infinity norm bounds for the inverse of Nekrasov matrices}

\author[pmf]{Ljiljana Cvetkovi\'c}
\author[kinkin]{Ping-Fan Dai}
\author[ftn]{Ksenija Doroslova\v cki}
\author[kin]{Yao-Tang Li}

\address[pmf]{Department of Mathematics and Informatics, Faculty of Science,
University of Novi Sad, Serbia}
\address[kinkin]{Department of Mathematics and Computer Science, Sanming University, Sanming, Fujian, 365004, PR China}
\address[ftn]{Faculty of Technical Sciences, University of Novi Sad, Serbia}
\address[kin]{School of Mathematics and Statistics, Yunnan University, Kunming, Yunnan, 650091, PR China}

\vspace{1cm}

\begin{abstract}
From the application point of view, it is important to have a good upper bound for the maximum norm of the inverse of a given matrix $A$.
In this paper we will give two simple and practical upper bounds
for the maximum norm of the inverse of a Nekrasov matrix. 
\end{abstract}

\end{frontmatter}

\maketitle

\section{Introduction}
\vspace{-3mm}
Throughout this paper, we will use the following notations:\\
$\mathbb C^n (\mathbb R^n) $ for the complex (real) $n$ dimensional  vector space,\\
$\mathbb C^{n,n}(\mathbb R^{n, n})$ for collection of all $n \times  n$ matrices with complex (real) entries,\\
$N:=\{ 1,2,\ldots,n\}$ for the set of indices,\\
$\displaystyle{r_i(A)=\sum_{j\in N \backslash \{i\}}^n |a_{ij}|},\;\;i\in N$ for the deleted  $i$-th  row sum.\\
Also, we define $h_i(A)$ recursively:
\begin{eqnarray*}
h_1(A) & := & \sum_{j \neq 1} |a_{1j}|, \\
h_i(A) & := & \sum_{j=1}^{i-1} |a_{ij}| \frac{h_j(A)}{|a_{jj}|} + \sum_{j=i+1}^{n} |a_{ij}|.
\end{eqnarray*}

We always denote by $A=D-L-U$ the standard splitting of a matrix $A$ into its diagonal $(D)$, strict lower
$(-L)$ and strict upper $(-U)$ triangular parts.

All classes of matrices, which are considered in this paper, are subclasses of $H$-matrices, more precisely nonsingular $H$-matrices. This class is
well-know from many applications, and it can be defined by:
\begin{definition}\label{H-mat}
A matrix $A=[a_{ij}]\in \mathbb C^{n,n}$is called an $H$-matrix if its comparison matrix $\langle A \rangle = [m_{ij}]$ defined by
$$\langle A\rangle=[m_{ij}]\in \mathbb C^{n, n};\;\;\;\;\;
m_{ij}= \left\{\begin{array}{ll}  \;\;\;|a_{ii}|, & i=j \\ -|a_{ij}|, & i\neq j \end{array}\right..$$

is an $M$-matrix, i.e., $\langle A \rangle ^{-1} \geq 0$.
\end{definition}

A well-known property of $H$-matrices is given by the following theorem (see \cite{BermanPlemmons}).
\begin{theorem}\label{Companion-mat}
Let $A=[a_{ij}]\in \mathbb C^{n, n}$ be a nonsingular $H$-matrix. Then
$$|A^{-1}|\leq \langle A  \rangle^{-1}.$$
\end{theorem}

The most important subclass of $H$-matrices is the class of strictly diagonally dominant ($SDD$) matrices, defined by:

\begin{definition}\label{SDD-mat}
 A matrix $A=[a_{ij}]\in \mathbb C^{n,n}$
is called SDD matrix if, for each $i \in N$, it holds that $$ |a_{ii}|> r_i(A). $$
\end{definition}

Finally, let us recall the famous result of Varah \cite{JM.Varah}, which gives us an upper bound for $SDD$ matrices:
 \be \| A^{-1} \|_{\infty} \leq \frac{1}{\displaystyle{\min_{i \in N}
(|a_{ii}| -r_i(A))}}. \label{Varah} \ee
Obviously, this upper bound works only for $SDD$ matrices, and even then it is not always good enough.
As we see in above written formula, the smaller $\displaystyle{\min_{i \in N} (|a_{ii}| -r_i(A))}$ is, the worse estimation will be. For that reason it can be
useful to obtain new upper bounds, for a wider class of matrices which, as we will see, will sometimes work better in the $SDD$ case, too.
The class which will be chosen for that purpose is the class of, so called, Nekrasov matrices.

\begin{definition}
A matrix $A=[a_{ij}]\in \mathbb C^{n,n}, n\geq 2$ is called Nekrasov matrix if, for each $i \in N$, it holds that
$$|a_{ii}|> h_i(A).$$
\end{definition}

It is known, see, for example \cite{CL}, \cite{Szulc}, that the class of Nekrasov matrices contains SDD class, while, on the other hand, it is a subclass of  $H$-matrices.

\section{New estimations}

We will present two possibilities to estimate the infinity norm of the inverse matrix of a Nekrasov matrix.
Although they look very similar, numerical examples will show that each of them can be better than the other one.\\

\begin{theorem}
Suppose  that $A=[a_{ij}]\in \mathbb C^{n, n}$ is a Nekrasov matrix. Then,
\be \displaystyle{\| A^{-1} \|_{\infty} \leq \frac{\displaystyle{\max_{i \in N} \frac{z_i(A)}{|a_{ii}|}}}{1 - \displaystyle{\max_{i\in N} \frac{h_i(A)}{|a_{ii}|}}}} \label{KIN ineq} \ee \noindent
and
\be\displaystyle{ \|A^{-1}\|_\infty\leq \frac{\displaystyle{\max_{i\in N}z_i(A)}}{\displaystyle{\min_{i\in N}(|a_{ii}|-h_i(A))}}\;,} \label{Nasa ineq} \ee \noindent
where
\be
z_1(A)  := 1, \;\; z_i(A) := \sum_{j=1}^{i-1} \frac{|a_{ij}|}{|a_{jj}|} z_j(A) +1 , \;\; i \in N \setminus \{1\}. \label{z}
\ee
\label{KIN}
\end{theorem}

In order to prove these estimations, we will start with the following lemma, proved by Robert in \cite{Robert}.

\begin{lm}\label{1}
Given any matrix $A=[a_{ij}]\in\mathbb C^{n, n},\;n\geq2$, with $a_{ii}\neq 0$ for all $i\in N$, then
\be h_i(A)=|a_{ii}|\big[ (|D|-|L|)^{-1} |U|e\big]_i,
\label{Shulc eq} \ee
\noindent where $e\in\mathbb C^{n}$ is the vector with all components equal to $1$.
\end{lm}

An immediate corollary of this Lemma is the following characterization of Nekrasov matrices, given by Szulc in \cite{Szulc}:  

\begin{theorem}\label{Nek-SDD}
$A$ matrix $A=[a_{ij}]\in\mathbb C^{n,n},\;n\geq2$ is a Nekrasov matrix if and only if
\be (|D|-|L|)^{-1} |U|e<e,
\label{Shulc1 eq} \ee \noindent i.e. if and only if
$E-(|D|-|L|)^{-1} |U|$ is an SDD matrix, where $E$ is the identity matrix..
\end{theorem}

{\em Proof of Theorem \ref{KIN}:} First of all, notice that elements of matrix $(|D|-|L|)^{-1} |U|$ are nonnegative,
because matrix $(|D|-|L|)$ is an $M$-matrix.  From inequality (\ref{Shulc1 eq}) we see that sum of all elements in each row is less then $1$, therefore we can conclude that all diagonal elements are also less than $1$.

Assume that matrix $A$ is a Nekrasov matrix. Then, from Theorem \ref{Nek-SDD}, we know that  $E-(|D|-|L|)^{-1} |U|$ is an SDD matrix, so we can apply Varah bound for estimation of the infinity norm of its inverse matrix:
$$\displaystyle{\| C^{-1}  \|_\infty\leq \max_{i\in N}\frac{1}{|c_{ii}|-r_i(C)}},$$
where $C:=E-(|D|-|L|)^{-1} |U|$.

According to the fact that all diagonal entries of matrix $(|D|-|L|)^{-1} |U|$ are less than $1$, we have:
$$|c_{ii}|= 1-\big[ (|D|-|L|)^{-1} |U|\big]_{ii},$$
$$r_i(C)=\sum_{j\neq i}^n \big[ (|D|-|L|)^{-1} |U|\big]_{ij},$$
and, therefore:
$$|c_{ii}|-r_i(C)= 1-\sum_{j=1}^n \big[ (|D|-|L|)^{-1} |U|\big]_{ij}=$$
$$=1-\sum_{j=1}^n \big[ (|D|-|L|)^{-1} |U|e\big]_i
=1-\displaystyle{\frac{h_i(A)}{|a_{ii}|}}.$$
Hence,
$$\displaystyle{\|C^{-1}  \|_\infty\leq \max_{i\in N}\frac{1}{1- \displaystyle{\frac{h_i(A)}{|a_{ii}|}}}} = \frac{1}{1 - \displaystyle{\max_{i\in N} \frac{h_i(A)}{|a_{ii}|}}}.$$
Now, in order to estimate $\|A^{-1}  \|_\infty$, it only remains to find
 a link between matrices $C^{-1}$ and $A^{-1}$.

Since 
$$C=(|D|-|L|)^{-1}\langle A\rangle,$$
it holds that 
$$\langle A\rangle=(|D|-|L|)C,$$
and, see  Theorem \ref{Companion-mat},

$$\|A^{-1}  \|_\infty\leq \|\langle A\rangle^{-1}  \|_\infty \leq\|C^{-1}\|_\infty \|(|D|-|L|)^{-1}  \|_\infty.$$

Finally, because  $|D|-|L|$ is an $M$-matrix,
$$\|(|D|-|L|)^{-1}\|_\infty =\|(|D|-|L|)^{-1}e\|_\infty, $$ 
and, if we denote by $y:=(|D|-|L|)^{-1}e, $ then $e=(|D|-|L|)y$, or, by components:
$$|a_{11}|y_1 = 1, \;\; |a_{ii}|y_i=1+\sum_{j=1}^{i-1}|a_{ij}| y_j ,\;\; i \in N \setminus \{1\}.$$
Since $|a_{ii}|y_i=z_i(A)$, see (\ref{z}), we get
 $$\|(|D|-|L|)^{-1}\|_\infty = \|y\|_{\infty} = \displaystyle{\max_{i\in N}\frac{z_i(A)}{|a_{ii}|}},$$ 
and the first upper bound (\ref{KIN}) is proved.

In order to prove the second upper bound, we will multiply matrix
$E-(|D|-|L|)^{-1} |U|$ with diagonal matrix $|D|$ from the left side, and denote: 
\be B:=|D|-|D|(|D|-|L|)^{-1} |U|.
\label{B eq} \ee \noindent
Clearly, this multiplication does not change SDD property, so matrix $B$ is also an SDD matrix, and we can, again, use Varah bound for estimation of the infinity norm of its inverse matrix:
$$\displaystyle{\| B^{-1}  \|_\infty\leq \max_{i\in N}\frac{1}{|b_{ii}|-r_i(B)}}.$$

Knowing that all diagonal entries of matrix $(|D|-|L|)^{-1} |U|$ are less than $1$, we have:
$$|b_{ii}|=|a_{ii}|-|a_{ii}|\big[ (|D|-|L|)^{-1} |U|\big]_{ii},$$
$$r_i(B)=\sum_{j\neq i}^n|a_{ii}|\big[ (|D|-|L|)^{-1} |U|\big]_{ij},$$
and, therefore:
$$|b_{ii}|-r_i(B)=|a_{ii}|-\sum_{j=1}^n|a_{ii}|\big[ (|D|-|L|)^{-1} |U|\big]_{ij}=$$
$$=|a_{ii}|-|a_{ii}|\big[ (|D|-|L|)^{-1} |U|e\big]_i
=|a_{ii}|-h_i(A).$$
Hence,
$$\displaystyle{\|B^{-1}  \|_\infty\leq \max_{i\in N}\frac{1}{|a_{ii}|-h_i(A)}} = \frac{1}{\displaystyle{\min_{i\in N}(|a_{ii}|-h_i(A))}}.$$
Now, in order to estimate $\|A^{-1}  \|_\infty$, it only remains to find
 a link between matrices $B^{-1}$ and $A^{-1}$.
Since 
$$B=|D|(|D|-|L|)^{-1}\langle A\rangle,$$
it holds that 
$$\langle A\rangle=(E-|L||D|^{-1})B,$$
and, see  Theorem \ref{Companion-mat},

$$\|A^{-1}  \|_\infty\leq \|\langle A\rangle^{-1}  \|_\infty \leq\|B^{-1}\|_\infty \|(E-|L||D|^{-1})^{-1}  \|_\infty.$$

Finally, in order to estimate $\|(E-|L||D|^{-1})^{-1}\|_\infty,$ we start with
$$\|(E-|L||D|^{-1})^{-1}\|_\infty =\|(E-|L||D|^{-1})^{-1}e\|_\infty, $$ 
which is true because  $E-|L||D|^{-1}$ is an $M$ matrix. 

If we denote $z(A):=(E-|L||D|^{-1})^{-1}e, $ then $e=(E-|L||D|^{-1})z(A)$, or, by components:

$$z_1(A) = 1, \;\; z_i(A)=1+\sum_{j=1}^{i-1}\frac{|a_{ij}|}{|a_{jj}|} z_j(A) ,\;\; i \in N \setminus \{1\},$$
which is, precisely, the definition of $z_i(A)$, given by (\ref{z}).

 Therefore, 
 $$\|(E-|L||D|^{-1})^{-1}\|_\infty = \|z(A)\|_{\infty} = \displaystyle{\max_{i\in N}z_i(A)},$$ 
and the second upper bound (\ref{Nasa ineq}) is also proved. $\Box$

Both estimations, given by  (\ref{KIN ineq}) and  (\ref{Nasa ineq}), can be used for bounding the infinity norm of the inverse of a Gudkov matrix, defined as a matrix for which there exists a permutation matrix $P$, such that $PAP^T$ is a Nekrasov matrix. However, before applying either  (\ref{KIN ineq}) or  (\ref{Nasa ineq}), one can {\em find} the permutation matrix $P$, in order to ensure the first row of $PAP^T$ to be strictly diagonally dominant.

\section{Numerical examples}

Keeping in mind that SDD matrices are subset of Nekrasov matrices, we will compare our two bounds with Varah's one, when it is applicable. For this purpose we will consider the following six matrices:

$A_1=\left[
\begin{array}{cccc}
-7  &  1 & -0.2 & 2\\
7  & 88  & 2   &  -3\\
2  & 0.5 & 13 & -2\\
0.5& 3 & 1    & 6
\end{array}
\right],$
$A_2=\left[
\begin{array}{cccc}
8   &  1 & -0.2 & 3.3\\
7   & 13  & 2   &  -3\\
-1.3& 6.7 & 13  & -2\\
0.5 &  3 & 1   & 6
\end{array}
\right],$

$A_3=\left[
\begin{array}{cccc}
21  &  -9.1 & -4.2 & -2.1\\
-0.7&  9.1  & -4.2 &-2.1\\
-0.7&  -0.7 &  4.9 &-2.1\\
-0.7& -0.7  & -0.7 &2.8
\end{array}
\right],$
$A_4=\left[
\begin{array}{cccc}
5  &  1 & 0.2 & 2  \\
1&  21  & 1 &-3\\
2&  0.5 &  6.4 &-2\\
0.5& -1  & 1 & 9
\end{array}
\right],$

$A_5=\left[
\begin{array}{ccc}
6&-3&-2\\
-1&11&-8\\
-7&-3&10
\end{array}
\right],$
$A_6=\left[
\begin{array}{cccc}
8 &  -0.5 & -0.5 & -0.5\\
-9 & 16 & -5 & -5\\
-6 &  - 4 &  15 & -3\\
-4.9 &  -0.9  & -0.9 & 6
\end{array}
\right].$

The upper bounds for their inverse matrices are given in the following table:

\begin{center}
\begin{tabular}{|c|c|c|c|c|c|} \hline
Matrix  & class & exact $\|A^{-1}\|_\infty$& Varah & (\ref{KIN ineq}) & (\ref{Nasa ineq}) \\\hline

$A_1$ & SDD & 0.1921& 0.6667 & 0.3805 & 0.5263  \\\hline
$A_2$ & SDD & 0.2390& 1 & 0.8848 &  0.6885   \\\hline
$A_3$ & SDD & 0.8759& 1.4286 &1.8076 & 0.9676    \\\hline
$A_4$ & SDD & 0.2707 & 0.5556 &0.6200 & 0.7937    \\\hline
$A_5$ & Nekrasov & 1.1519& - & 1.4909& 2.4848  \\\hline
$A_6$ & Nekrasov & 0.4474& - & 1.1557& 0.5702  \\\hline
\end{tabular}
\end{center}

Obviously, each estimation (\ref{KIN ineq}) or (\ref{Nasa ineq}) can work better than the other one. So, in general case, for Nekrasov matrices, one can take the smallest estimation of these two.

For SDD case, it can be, also, useful to invest more time and calculations in obtaining these new estimations, since they can be significantly better (for example, see matrix $A_2$).

\section{Acknowledgments}
This work is partly supported by the Ministry of Science and Environmental Protection of Serbia (174019),  by the Provincial Secretariat of Science and Technological Development of Vojvodina, Serbia (2002), and   by the Natural Science Foundation of the Education Department of Fujian province, China (JA11256, JB09228).

\end{document}